\documentclass[ijaa]{informs4} 
\DoubleSpacedXI 
\usepackage{eqndefns-left} 
\RequirePackage{tgtermes}
\RequirePackage{newtxtext}
\RequirePackage{newtxmath}
\RequirePackage{bm}
\RequirePackage{endnotes}
\OneAndAHalfSpacedXII
\usepackage{graphicx}
 \usepackage{float}
\usepackage{amsmath}
\usepackage{caption}
\usepackage{subcaption}

\newcommand{\E}{\mathbb{E}}
\newcommand{\setO}{\mathbf{O}}
\newcommand{\N}{\mathbf{N}}
\newcommand{\setM}{\mathbf{M}} 
\newcommand{\setR}{\mathbf{R}} 
\newcommand{\setI}{\mathbf{I}} 
\newcommand{\setD}{\mathbf{D}} 

\newcommand{\setS}{\mathbf{S}}
\newcommand{\setE}{\mathbf{E}}
\usepackage{makecell}
\usepackage{endnotes}
\usepackage{booktabs}
\usepackage{bbm}

\usepackage{hyperref}
\usepackage{natbib}
 \bibpunct[, ]{(}{)}{,}{a}{}{,}%
 \def\bibfont{\small}%
\TheoremsNumberedThrough     
\ECRepeatTheorems

\EquationsNumberedThrough    

\usepackage{makecell, multirow}
\usepackage{optidef}

\usepackage[ruled,vlined, noend,linesnumbered]{algorithm2e}

\begin{document}


\RUNAUTHOR{Shen et al.} 

\RUNTITLE{Data-driven assortment and inventory at JD.com}

\TITLE{JD.com Improves Fulfillment Efficiency with Data-driven Integrated Assortment Planning and Inventory Allocation}
\ARTICLEAUTHORS{%
\AUTHOR{Zuo-Jun Max Shen\textsuperscript{1,2}, Shuo Sun\textsuperscript{1}}
\AFF{\textsuperscript{1}Department of Industrial Engineering and Operations Research, UC Berkeley, Berkeley, CA 94720 \\
\AFF{\textsuperscript{2}Faculty of Engineering and Faculty of Business and Economics,
The University of Hong Kong, Hong Kong, China}
\EMAIL{maxshen@hku.hk}, \EMAIL{shuo\_sun@berkeley.edu}}
\AUTHOR{Yongzhi Qi\thanks{Corresponding Author}, Hao Hu, Ningxuan Kang, Jianshen Zhang, Xin Wang, Xiaoming Lin}
\AFF{Supply chain tech team, JD.com, Beijing 101111, China \\
\EMAIL{qiyongzhi1@jd.com}, \EMAIL{huhao@jd.com}, \EMAIL{kangningxuan@jd.com},
\EMAIL{zhangjianshen@jd.com},
\EMAIL{wangxin52@jd.com},
\EMAIL{linxiaoming7@jd.com}
}
} 
\ABSTRACT{This paper presents data-driven approaches for integrated assortment planning and inventory allocation that significantly improve fulfillment efficiency at JD.com, a leading E-commerce company. JD.com uses a two-level distribution network that includes regional distribution centers (RDCs) and front distribution centers (FDCs). Selecting products to stock at FDCs and then optimizing daily inventory allocation from RDCs to FDCs is critical to improving fulfillment efficiency, which is crucial for enhancing customer experiences. For assortment planning, we propose efficient algorithms to maximize the number of orders that can be fulfilled by FDCs (local fulfillment). For inventory allocation, we develop a novel end-to-end algorithm that integrates forecasting, optimization, and simulation to minimize lost sales and inventory transfer costs. Numerical experiments demonstrate that our methods outperform existing approaches, increasing local order fulfillment rates by 0.54\% and our inventory allocation algorithm increases FDC demand satisfaction rates by 1.05\%. Considering the high-volume operations of JD.com, with millions of weekly orders per region, these improvements yield substantial benefits beyond the company's established supply chain system. Implementation across JD.com's network has reduced costs, improved stock availability, and increased local order fulfillment rates for millions of orders annually.}
\maketitle
\section*{Introduction}
JD.com, one of the largest E-commerce companies in the world, has significantly enhanced its fulfillment efficiency by implementing data-driven assortment planning and inventory allocation strategies. Leveraging artificial intelligence and advanced analytics, JD.com has optimized its supply chain and order fulfillment processes, resulting in increased customer satisfaction and business growth.

\begin{figure}[htbp]
  \setlength{\abovecaptionskip}{0pt}
  \setlength{\belowcaptionskip}{0pt} 
  \begin{minipage}[b]{0.53\textwidth}
    \centering
    \caption{The Graphic Illustrates a Two-level Inventory Network}
    \includegraphics[width=\linewidth]{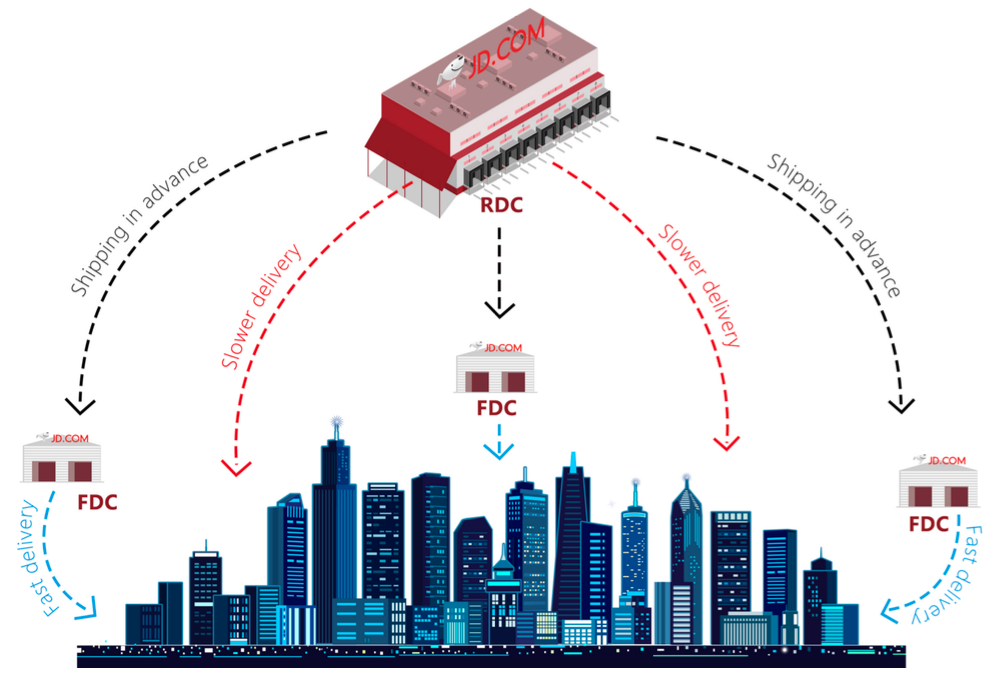}
    \label{fig:inv}
  \end{minipage}  \hspace{40pt} 
  \begin{minipage}[b]{0.33\textwidth}
    \centering
    \caption{The Graphic Visualizes the Decision Flow}
    \includegraphics[width=\linewidth]{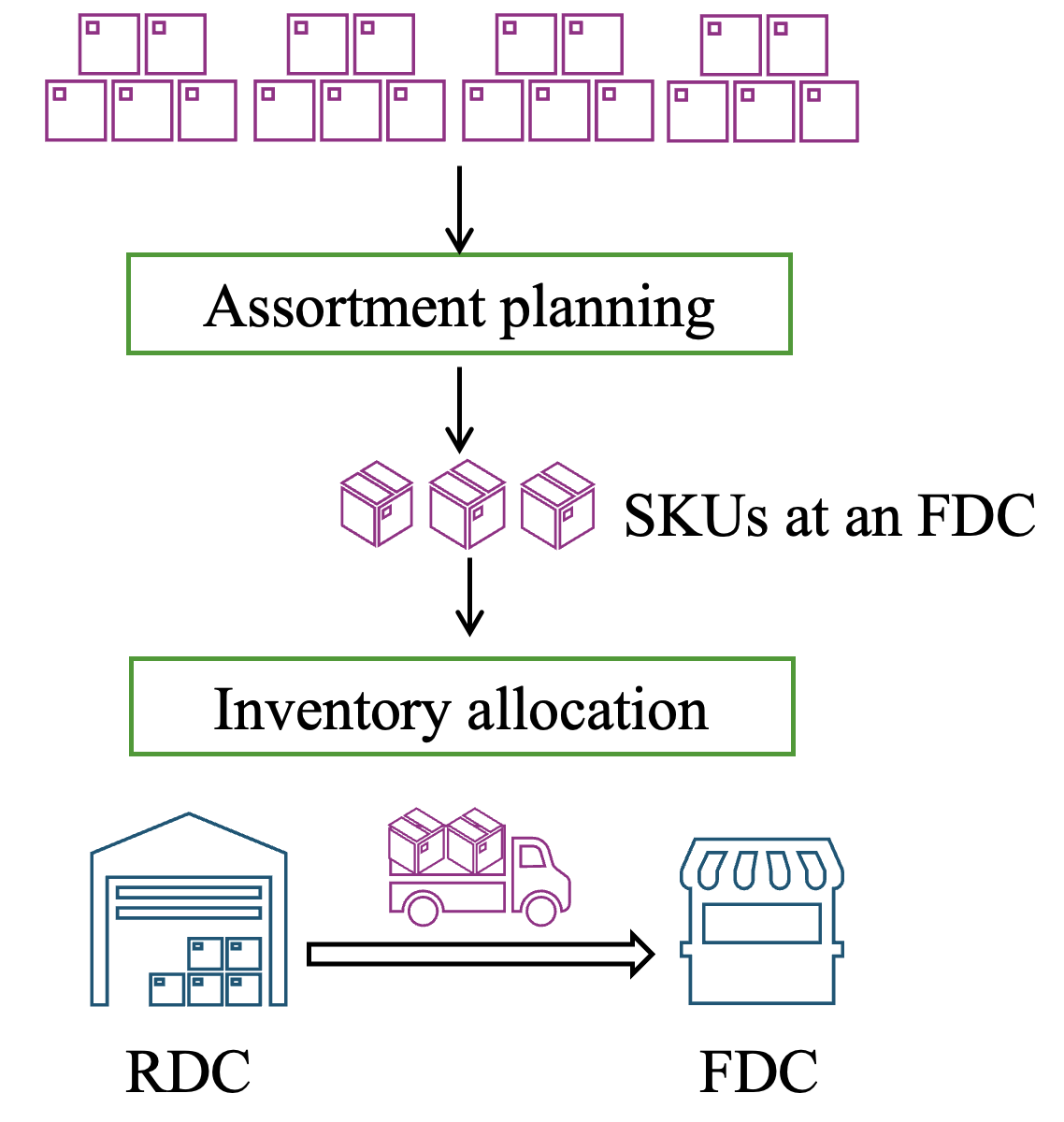}
    \label{fig:alg}
  \end{minipage}
\end{figure}

Fulfillment efficiency is crucial for ensuring customer satisfaction in the E-commerce industry. JD.com excels in this area due to its extensive logistics infrastructure and commitment to delivering 90\% of all orders within 24 hours to over 588 million active users \citep{JD2023AnnualReport}. This rapid delivery promise has significantly enhanced customer satisfaction and distinguished JD.com from its competitors. 
JD.com employs a two-level delivery network comprising Regional Distribution Centers (RDCs) and Front Distribution Centers (FDCs). The eight RDCs across China are large warehouses stocking all stock-keeping units (SKUs) sold online, and serving broad regions. Closer to customers, FDCs provide rapid delivery but have smaller capacities and rely on daily shipments from RDCs. This two-level inventory network is illustrated in Figure \ref{fig:inv}. When an order is placed, it is fulfilled from the corresponding FDC if sufficient inventory is available. If the FDC is out of stock, the RDC fulfills the order, incurring extra costs due to longer distances and potential lost sales from extended delivery times. 
Given the importance of FDCs for fast delivery and their capacity constraints, the selection of products to store at FDCs and the inventory management at FDCs and RDCs are critical for maintaining fulfillment efficiency. Our goal is to maximize fulfillment efficiency, minimize sales losses at FDCs and RDCs, and reduce inventory transfer costs. 

The problem is challenging for the following five reasons: (1) A single RDC serves multiple FDCs and has its own demand requirements, necessitating intricate coordination of its own demand fulfillment and inventory transfers to various FDCs. This also involves various business constraints, such as capacity constraints at FDCs, logistics throughput, and sales plans. (2) JD.com manages a vast number of SKUs across different product categories, requiring the algorithm to be scalable. (3) The demand at FDCs is relatively sparse and variable, complicating demand prediction. (4) Assortment and inventory decisions must account for future uncertainties to align with long-term sales plans and mitigate risks from unforeseen events. (5) The interpretability of the prediction algorithm is crucial for coordinating different sections of the supply chain. To address these challenges, we propose an integrated optimization model aimed at maximizing the number of orders fulfilled by FDCs while considering various business constraints. Due to the complexity of the problem, even a static version with deterministic demand is NP-complete. 
Additionally, the problem size is substantial since each RDC manages over 10 million SKUs. Consequently, we decompose the problem into two subproblems: the {\it assortment planning} problem and the {\it inventory allocation} problem. The assortment planning problem determines which SKUs to stock at FDCs, which in turn affects the warehouse capacity usage and transportation needs in the inventory problem. The assortment decisions also impact the achievable service level at FDCs. Furthermore, the quality of the assortment solution has a direct impact on the solution quality of the inventory problem. We visualize the decision flow in Figure \ref{fig:alg}.

In the {\it assortment planning} subproblem, we select which SKUs to stock at FDCs over a specific period (e.g., one week) to maximize the number of fulfilled orders at FDCs under capacity constraints. This subproblem presents two key challenges: predicting future order structures is difficult, and even with precise future order information, the resulting integer programming problem is NP-hard and even challenging to approximate. This optimization problem can be considered as a special case of maximizing a supermodular function under capacity constraints, which is generally hard to solve. To address these challenges, we propose three heuristic algorithms: \emph{ML-Top-K}, \emph{Reverse-Exclude}, and \emph{Hybrid-Selection}. The ML-Top-K algorithm utilizes machine learning (ML) techniques to predict the number of orders containing each product, and then ranks the products in descending order based on these predicted order counts, selecting the top $K$ products from this ranking. The Reverse-Exclude algorithm iteratively eliminates products from the consideration set based on their historical order frequency. The Hybrid-Selection algorithm takes advantage of the output of the previous two heuristics. We compare local order fulfillment rates of the algorithms in numerical experiments, which refer to the percentage of customer orders that can be satisfied directly from FDC stock without order fulfillments by RDCs.
Numerical experiments demonstrate that all three methods achieve higher local order fulfillment rates compared with previous methods. The assortment optimization controls the number of SKUs, ensuring that the downstream inventory problem naturally satisfies some constraints.

After determining which SKUs to store in the FDCs, we address the second subproblem, the {\it inventory allocation} problem, where the daily transfer amount of on hand inventory from RDCs to FDCs is determined to minimize total sales losses and inventory transfer costs. We model this problem as a periodic review inventory problem. 
This problem relates to the two-echelon inventory allocation and rationing problem, where a supplier distributes available on hand inventory among retailers \citep{snyder2019fundamentals}. However, our problem differs because the supplier (RDC) also needs to allocate sufficient inventory to itself. This subproblem is challenging due to demand uncertainty, the interplay between FDCs and RDCs, the interdependence of inventory decisions across multiple time periods, and various capacity and business constraints. To address these challenges, we propose a novel multitask end-to-end learning framework. This framework outperforms previous algorithms due to its multitask and self-learning capabilities, as well as its interpretability. We provide further details in the {\it Contributions} subsection and in the {\it End-to-End Algorithm for Inventory Allocation} section.

\subsection*{Contributions}
\label{sec:contribution}
Our work makes the following key contributions:

1. We propose an integrated framework for joint assortment planning and inventory allocation, aimed at improving the number of orders fulfilled by FDC and minimizing lost sales due to stockouts at both FDCs and RDCs. 

2. For the assortment planning problem, we introduce three novel algorithms that maximize FDC order fulfillment, addressing the challenges of order structure prediction and optimization complexity. The assortment planning problem, which aims to maximize local fulfillment for FDCs, has been relatively understudied in the previous literature. We test the performance of our algorithms using large-scale real data.  

3. For the inventory allocation problem, we develop a novel end-to-end forecasting and allocation algorithm. This algorithm combines forecasting with a simulation system to learn target inventory (TI) and safety stock (SS) levels. Unlike previous algorithms that rely solely on historical data for optimal inventory labeling, our approach integrates forecasting and simulation to minimize total loss \citep{qi2023practical}. This learning framework offers several advantages.
\begin{itemize}
    \item \emph{Multitask learning.} Our framework enables the network to simultaneously learn inventory allocation decisions across FDCs and RDC, and demand prediction to minimize the total cost and demand prediction error. This facilitates information sharing and knowledge transfer between tasks and thus improves the efficiency and accuracy of inventory management. 
    \item \emph{Self-learning.} We introduce an advanced self-learning mechanism that eliminates the need for complex operations research algorithms previously used for data labeling. This significantly enhances the efficiency and accuracy of sample processing in inventory management.
    \item \emph{Interpretability.} The modular structure of our model enhances interpretability compared with traditional "black-box" deep learning methods, providing better insights into the inventory management process and facilitating more informed decision making. We provide more details in the {\it End-to-End Algorithm for Inventory Allocation} section.
\end{itemize}

\subsection*{Numerical Experiment, Implementation and Impact}
 We have developed, tested, and implemented a decision support system utilizing these algorithms across JD.com's inventory network in China. The numerical experiments show that joint assortment and inventory optimization can efficiently satisfy customers' demands, which significantly improves the local order fulfillment rate at the FDCs.
 Our decision system has substantially enhanced operational efficiency and customer experiences. Key benefits include reduced inventory holding and transfer costs, improved stock availability, and increased local order fulfillment rates. This improvement in the local order fulfillment rate is particularly beneficial for JD.com's speedy `211 program', which guarantees same-day delivery for orders placed before 11 am and next-day delivery before 3 pm for orders submitted before 11 pm. We provide details in the {\it Implementation and Impact} subsection.  
 
\section*{Related Work}
Our work relates to (warehouse) assortment optimization, inventory allocation in a two-echelon setting, and machine learning in inventory management.

Assortment planning has been extensively studied in the context of revenue management, where it refers to determining which products to display to customers in retail settings to maximize total revenue  \citep{ryzin1999relationship}. This stream of research typically assumes that customers select products based on a discrete choice model. Various choice models and constraints have been explored in retail assortment optimization; \cite{kok2015assortment} and \cite{gallego2019assortment} provide comprehensive reviews. We refer to this body of work as \textit{retail assortment planning}. 
In contrast, our problem falls under \textit{warehouse assortment planning}, which differs from retail assortment planning in both demand modeling and optimization objectives: First, demand modeling in warehouse assortment planning cannot be effectively captured by discrete choice models. Unlike retail settings, where products are displayed to customers simultaneously, warehouse assortment involves selecting products to stock in FDCs. So customers would not see all the products in the assortment at the same time. Moreover, customers may purchase multiple products and units in one order, with complex dependencies between products. Although multiple-purchase choice models have been studied (e.g., \citet{tulabandhula2023multi, bai2024assortment, jasin2024assortment}), the estimation and optimization become intractable as the number of products in an order grows. Furthermore, demand in warehouses is highly dynamic, influenced by promotions and marketing events, making discrete choice models unsuitable. Instead, we employ an aggregate demand function model. Second, our optimization objective is to maximize fulfillment efficiency rather than total expected revenue, which is the traditional focus in retail assortment planning.

Research on warehouse assortment planning is relatively limited. \cite{catalan2012assortment} first studied assortment selection at distribution centers with the goal of minimizing split orders when demand is known. They showed the problem is NP-hard and proposed heuristics to solve it. \cite{li2024multi} extended this work to multiple-warehouse assortment selection, again focusing on minimizing order splits under known order distributions. \cite{shen2024near} used distributionally robust optimization to address multiple-warehouse assortment planning, considering demand ambiguity and proposing a heuristic based on item popularity and warehouse costs. For single-warehouse assortment selection, \cite{wu2019assortment} aimed to maximize the order fulfillment rate under cardinality constraints, using exponential smoothing and graph-theoretic heuristics to predict demand and formulating a robust optimization model. However, in E-commerce settings like JD.com, accurately predicting the demand of different orders at FDCs is impractical, and solving large-scale mixed-integer programs becomes computationally intractable. 
\cite{li2022shall} is the most closely related work, addressing single-warehouse assortment selection under cardinality constraints to minimize order fulfillment costs. When order types and their demand distributions are known, our problem reduces to a special case of theirs. However, accurate estimates of order types and demands are infeasible in large-scale E-commerce settings. Although their heuristic, marginal choice indexing, is equivalent to our ML-top-K algorithm when product demand is known, we use deep learning to predict demand and provide a new theoretical guarantee. Additionally, we propose novel heuristics, including Reverse-Exclude and a hybrid algorithm, and test our approaches on extensive real-world data from JD.com.
 
Unlike these works, most of which focus on product selection, we consider the downstream inventory allocation decisions. Although recent research has begun addressing joint assortment and downstream decisions such as inventory and matching \citep{bai2022coordinated, huang2024basic, epstein2024optimizing, sun2025unifiedalgorithmicframeworkdynamic}, these studies remain within the field of retail assortment optimization.  More specifically, these works focus on maximizing the total expected revenue and use discrete choice models to model demand while we aim to maximize fulfillment efficiency and use aggregate demand functions. Instead, we show the success of assortment optimization in improving inventory management and fulfillment efficiency.  

Our paper is also related to inventory allocation in a two-echelon setting, where a supplier distributes available on hand inventory among retailers \citep{snyder2019fundamentals}. An RDC differs from classic inventory allocation problems in that it must keep inventory for itself to fulfill orders that FDCs cannot fulfill. Two-echelon inventory allocation problems have been widely studied after the seminal work of \cite{clark1960optimal}. Under known demand distributions, various solution algorithms have been proposed for different objectives and constraints \citep{Sherbrooke1968,tan1974optimal,caglar2004two,huang2005two}. In cases of unknown demand distributions, \cite{nambiar2021dynamic} developed a Lagrangian heuristic to transform the problem into a single-echelon inventory management problem under demand learning using a Bayesian demand model and ARMA/ARIMA model. Motivated by the applications in the fashion industry, \cite{fisher1996reducing, fisher2000accurate}, and \cite{gallien2015initial} have focused on two-period inventory allocation problems. However, these approaches are unsuitable in settings with highly variable demand where traditional demand update methods fail. Also, we consider a multiperiod inventory allocation problem. Our problem also intersects with joint inventory and fulfillment optimization \citep{jasin2015lp,amil2022multi, ma2023order}. These studies assume known future order distributions, which is challenging in our E-commerce setting due to the unpredictability of the order structure. 

Finally, recent advances in machine learning for inventory management have applied end-to-end (E2E) learning and reinforcement learning (RL) techniques. In RL, \cite{oroojlooyjadid2022deep} developed a deep reinforcement learning (DRL) algorithm for the beer game, while \cite{gijsbrechts2022can} applied the Asynchronous Advantage Actor Critic (A3C) algorithm to dual-sourcing multiple-echelon problems. \cite{liu2022multi} explored multiple-agent DRL for decentralized multiple-echelon inventory management. In E2E learning, \cite{qi2023practical} proposed a single-warehouse replenishment framework using dynamic programming to label optimal replenishment quantities and use deep learning to determine the optimal replenishment decisions. Although this method has been successfully implemented at JD.com, it cannot be applied to our two-echelon inventory allocation problem due to the intractability of computing the labels (i.e., the optimal allocation decisions), which arises from two main challenges: (1) The optimal allocation for each period cannot be computed via dynamic programming because of the large state space, corresponding to the inventory levels at the RDC and multiple FDCs. (2) Solving the dynamic program requires addressing multiple mixed-integer programs, which lack closed-form solutions. Furthermore, JD.com has observed that the original E2E algorithm performs poorly under fluctuating demand. Moreover, our algorithm offers improved interpretability, which is critical for practical applications in dynamic E-commerce settings.

\section*{Problem Formulation}
In this section, we formally define the assortment and inventory optimization problem. We start with the assortment subproblem, and then define the inventory allocation problem. 

\subsection*{Assortment Problem}
Let $\N$ denote the set of $N$ distinct SKUs, indexed $1,2,\dots, N$. In the assortment planning problem, we select a subset of SKUs $\setS \subseteq \N$ to store at the FDCs. The set $\setS$ includes at most $K$ SKUs, which reflects a practical business requirement. Our goal is to maximize the number of orders that can be fulfilled at FDCs during the selling horizon. 

Let $\setO$ denote the set of order types. Each order type $o \in \setO$ represents a unique combination of SKUs that are ordered together. Let $D_o$ denote the number of orders corresponding to order type $o$ in a given future selling horizon. Let $S_o$ denote the set of SKUs that is included in order type $o$. Then the problem can be formulated as the following 0-1 integer program:

\begin{subequations}
\label{model:selection}
\begin{align}
\max_{x,y}\;& \sum_{o\in\setO} D_o\,y_o \label{model:obj}\\
\text{s.t.}\;& \sum_{i\in\N} x_i \;\le\; K, \label{model:card}\\
& y_o \;\le\; x_i,
  &&\forall\,o\in\setO,\;i\in S_o, \label{model:link}\\
& x_i \in \{0,1\},\; y_o \in \{0,1\},
  &&\forall\,i\in\N,\;o\in\setO. \label{model:bin}
\end{align}
\end{subequations}

Here, $x_i$ is a binary variable that indicates whether the product $i$ is included in the assortment, and $y_o$ is a binary variable that indicates whether the FDC can fulfill the type of order $o$. The first constraint requires that at most $K$ unique SKUs can be included in the assortment. The second requires that one order can be fulfilled by the FDC only if all its products are selected in the assortment.

Solving this problem is challenging for two reasons. First, accurately forecasting future order structures is difficult due to high-dimensional sparse data, nonstationary demand patterns, and complex product interdependencies inherent in E-commerce. Second, the optimization problem, given the order structure, is NP-hard, making it computationally intractable to find exact solutions at the scale of typical E-commerce operations.

\subsection*{Inventory Allocation}
In the inventory allocation algorithm, we consider a multiperiod inventory allocation problem with stochastic demand and no lead time, over a finite horizon of discrete periods $1, \dots, T$. 
In our two-level inventory allocation network, each FDC can be replenished by only one RDC. In each period, JD.com determines the replenishment quantities at each FDC. When the selling horizon starts, orders arrive sequentially at both the FDCs and the RDC, and we must decide whether to fulfill them using the FDC or the RDC. Each FDC has a service region, and orders arriving in this region can only be fulfilled by the FDC and its corresponding RDC. Since inventory and fulfillment decisions are independent across different RDC regions, we only need to consider one RDC and its FDCs.

We model the problem as a periodic review inventory problem over $T$ periods. Let $\mathbf{J}$ denote the set of FDCs. Subscript 0 refers to variables and parameters associated with the RDC. Let $D_{ij}^t$ denote the demand for product $i$ at FDC $j$ in period $t$, and let $\mathbf{D}^t = (D_{ij}^t)_{i \in \N, j \in \mathbf{J}\cup\{0\}}$. The lead time for inventory transfers is denoted by $l$. Consider any product $i$. Let $I_{ij}^t$ denote the on-hand inventory of product $i$ at FDC $j$ at the start of period $t$, and let $I_{i0}^t$ be the inventory at the RDC after replenishment at the beginning of period $t$. The RDC replenishment decision is outside the scope of this work and is considered exogenous. In each period $t$, FDC $j$ receives inventory $u_{ij}^{t-l}$, which is the amount transferred $l$ periods ago, and the RDC sends out new inventory $u_{ij}^t$ to each FDC. Orders then arrive at both the FDCs and RDC. For an order at FDC $j$, only inventory from FDC $j$ and the RDC can be used to fulfill it. Moreover, no backorders are allowed.

Let the state at the start of period $t$ be \(
    S_t \;=\; \bigl(\mathbf I^t,\; \mathbf U^t\bigr)\),
where $\mathbf U^{\,t}=\{u_{ij}^{k}\}_{i\in\N,\;j\in\mathbf J, \;k\in \{t-1,\dots t-l\}}$ stores the most recent $l$ shipments that have not yet reached the FDCs. In period $t$, we need to determine a nonnegative transfer vector $\mathbf{u}^t$, which will arrive at the FDCs at the beginning of period $t+l$. 

Our goal is to minimize the total cost in $T$ periods, which includes lost-sales costs at both the RDC and FDCs, inventory-transfer costs, and the cost of serving any unfilled FDC orders from the RDC. We first consider period $t$. Let $c$ denote the unit cost of fulfilling FDC demand from the RDC, and let $s$ denote the unit cost of demand left unsatisfied by both the FDC and the RDC. Consider any product $i\in \N$. Let $r_{ij}$ denote the per-unit cost of inventory transfer between FDC $j$ and RDC. We define the decision variables as follows: $x_{ij}^t$ denotes the number of product $i$ fulfilled by FDC $j$, and let $y_{ij}^t$ be the number of units of product $i$ from demand at FDC $j$ that is fulfilled by the RDC, $y_{i0}^t$ denotes the fulfilled demand of RDC. $z_{ij}^t$ denotes the lost sales of product $i$ from demand at FDC $j$, and $z_{i0}^t$ denotes lost sales at RDC. The total cost in period $t$ is computed via the following linear program:
\begin{align}
C^t(\mathbf{D}^t,\mathbf{u}^{t-l},\mathbf{u}^t, \mathbf{I}^t) = &\min_{\mathbf{x}, \mathbf{y}, \mathbf{z}} c \sum_{j\in\mathbf{J}}y_{ij}^t + s\sum_{i\in\mathbf{N}} \sum_{j\in\mathbf{J} \cup\{0\}}z_{ij}^t \label{eq:objective}  \\
\text{s.t.} \quad &x_{ij}^t\leq I_{ij}^t + u_{ij}^{t-l} , \quad \forall j\in\mathbf{J},  i\in\N \label{eq:fdc_inventory} \\
&  \sum_{j\in \mathbf{J}}y_{ij}^t +y_{i0}^t+\sum_{j\in \mathbf{J}} u_{ij}^t\leq I_{i0}^t, \quad \forall i\in\N \label{eq:rdc_inventory} \\
& x_{ij}^t + y_{ij}^t + z_{ij}^t = D_{ij}^t, \quad \forall j\in\mathbf{J}, i\in \N \label{eq:demand_balance} \\
& y_{i0}^t +z_{i0}^t =D_{i0}^t, \quad \forall i\in \N  \\
& x_{ij}^t, y_{ij}^t, z_{ij}^t\geq 0, \quad \forall j\in\mathbf{J}, i\in\N \label{eq:non_negativity} \\
& \text{Other business constraints}. \nonumber
\end{align}
The objective accounts for lost sales at both an RDC and FDCs, and the added cost of fulfilling FDC demand from the RDC. Constraint \eqref{eq:fdc_inventory} ensures that the number of products used to fulfill orders at FDCs does not exceed the available inventory plus the transferred amount. Constraint \eqref{eq:rdc_inventory} limits the fulfillment from the RDC and the transfer to FDCs by its available inventory. Constraint \eqref{eq:demand_balance} requires that the sum of fulfilled demand and lost sales equals the total demand at FDCs or the RDC. 

Next, we define the state transitions. The transfer buffer updates as $
    \mathbf U^{t+1}
    \;=\;
    \bigl(\mathbf u^{t},\mathbf u^{\,t-1},\dots,\mathbf u^{\,t-l+1}\bigr)$. 
For $\mathbf I_t$, the optimal solution $(\mathbf{x},\mathbf{y}, \mathbf{z})$, the on-hand inventory level at the start of period $t+1$ is given by $I_{ij}^{t+1}=\max\{I_{ij}^t +u_{ij}^{t-l}- x_{ij}^t,0 \}$ for $j\in \mathbf{J}$ and $i\in \mathbf{N}$, $I_{i0}^{t+1}= \max\{I_{i0}^t-\sum_{j\in \mathbf{J}\cup\{0\}} y_{ij}^t,0 \}+ R_i^{t}$, where $R_i^t$ denotes the uncertain replenishment from suppliers. Let $\mathbf R^t=(R_i^t)_{i\in \mathbf{N}}$. With the terminal value $V_{T+1}(\cdot)=0$, we define the optimal cost-to-go function $V_t\bigl(\mathbf I^{t},\mathbf U^t\bigr)$ for each period $t=T,T-1,\dots,1$, which represents the minimum expected total cost from period $t$ to the end of the planning horizon. 
\begin{equation}\label{eq:bellman}
    V_t\bigl(\mathbf I^{t},\mathbf U^t\bigr)
    \;=\;
    \min_{\mathbf u^{t}\ge 0}\;
    \E_{\mathbf D^{t},\mathbf R^{t}}\!\Bigl[C^t(
      \mathbf{D}^t,\mathbf{u}^{t-l},\mathbf{u}^t, \mathbf{I}^t)
         +\sum_{i\in \N,j\in \mathbf{J}}r_{ij}u_{ij}^t+ \sum_{i\in \mathbf{N}, j\in \mathbf{J}} V_{t+1}\bigl(\mathbf I^{\,t+1},\mathbf U^{t+1}\bigr)
    \Bigr].
\end{equation}

This problem is challenging for several reasons. First, accurately forecasting demand $D_{ij}^t$ across multiple FDCs and time periods is difficult due to high variability and sparsity. This challenge is exacerbated by the large number of SKUs and locations. Moreover, the replenishment quantity at the RDC is uncertain, requiring real-time decisions based on current system states. With perfect demand and replenishment information, although solving the LP for $C^t(\mathbf{D}^t,\mathbf{u}^{t-l},\mathbf{u}^t, \mathbf{I}^t)$ is straightforward via a greedy allocation rule that prioritizes using local inventory, but solving the dynamic program to determine the optimal $\mathbf{u}$ across many SKUs is computationally intensive due to the large state space. In practice, JD.com assumes deterministic forecasts for demand and replenishment, reformulates the dynamic program as a linear program, and still faces solution times on the order of hours. Such delays are problematic in the fast-paced E-commerce setting, where continuous order arrivals can quickly shift the inventory state, rendering previously computed solutions suboptimal. These limitations illustrate the challenges of applying traditional optimization techniques to large-scale, dynamic E-commerce inventory allocation. To address these issues and enable fast, adaptive decision making, we propose an end-to-end learning algorithm that integrates demand forecasting, inventory optimization, and fulfillment decisions into a unified framework. This approach is designed to rapidly adapt to evolving conditions while managing the scale and complexity inherent in modern E-commerce inventory systems.

\section*{Algorithms for Assortment Optimization}
In this section, we introduce three heuristics to solve the assortment optimization problem. Given the inherent difficulty of predicting order structures, we can only forecast the \emph{order frequency} for product $i$, defined as the number of distinct orders containing product $i$. 
The first algorithm uses machine learning to predict the order frequency of each product and selects the top $K$ products. We call this heuristic \emph{ML-Top-K}. We demonstrate that this greedy approach performs well when the ratio of single-item orders is high; see counterexamples in \cite{li2022shall}. However, this algorithm can perform arbitrarily poorly when the percentage of multiple-item orders is high. The key issue is that products with high order quantities may frequently appear in orders alongside products with low order quantities; therefore, simply selecting the top $K$ products cannot guarantee a high order fulfillment ratio. 

Motivated by this observation, we propose a second heuristic called \emph{Reverse-Exclude}. This algorithm progressively eliminates products from the ground set, which is the complete set of SKUs potentially selected at the FDC, starting with those that appear least frequently in orders. When a product is removed from the consideration set, we also eliminate all historical orders containing that product. This approach ensures that all remaining orders can be satisfied by selecting from the consideration set. Finally, we propose a hybrid method that takes advantage of the assortments from the previous two heuristics.
\subsection*{ML-Top-K Algorithm}
This algorithm first predicts product order quantities and then selects the top K products with the highest predicted quantities. We begin by describing the ML algorithm and then establish theoretical results under the assumption of accurate predictions. 

The ML-Top-K algorithm's input includes historical order quantities at various levels: SKU-specific, SKU category, brand, and nationwide SKU-level, along with promotional marketing activity factors. The model's core structure consists of seasonal and trend blocks to learn baseline patterns. Outputs from these blocks are processed through two parallel paths: a Temporal Convolutional Network (TCN) for deep feature extraction and a Multilayer Perceptron (MLP) for multiple-level feature extraction. The TCN output features are combined with the promotional predictions from the MLP through element-wise multiplication. This combined output is then processed through another MLP to generate the final order quantity forecast.
Based on this forecast, we select the top $K$ products with the highest predicted order quantities. Next, we show the theoretical guarantee of this greedy algorithm when the order frequency can be accurately forecasted. 
\begin{theorem}
Let $\alpha$ denote the minimum ratio of single-item orders to total orders including the item, among the top $K$ products when ranked by total order frequency. Assuming accurate order frequency prediction, the Top-K algorithm provides an $\alpha$-approximate solution. 
\label{thm:1}
\end{theorem}
\proof{Proof of Theorem \ref{thm:1}.}
Let $ALG$ and $OPT$ denote both the assortments selected by the greedy algorithm and the optimal solution, respectively, as well as the number of orders they can satisfy. Let $f_i$ denote the total number of orders containing product $i$. 

The proof follows in three steps. First, by definition of $\alpha$, for the assortment selected by the greedy algorithm:
\[
ALG \geq \alpha \sum_{i \in ALG} f_i,
\]
where $\alpha$ represents the minimum ratio of single-item orders to total orders among selected products.

Second, because the greedy algorithm selects the top $K$ products with the highest order frequencies:
\[
\sum_{i \in ALG} f_i \geq \sum_{i \in OPT} f_i.
\]

Third, for any assortment, including $OPT$, the total number of satisfied orders cannot exceed the sum of order frequencies of its products:
\[
\sum_{i \in OPT} f_i \geq OPT.
\]

Combining these inequalities:
\[
ALG \geq \alpha \sum_{i \in ALG} f_i \geq \alpha \sum_{i \in OPT} f_i \geq \alpha OPT.
\]

Therefore, the top-K algorithm achieves an approximation ratio of $\alpha$.
\Halmos
\endproof
However, when $\alpha$ is small, we have counterexamples showing that the greedy algorithm performs poorly. This shortcoming of greedy approaches has also been noted by \citep{li2022shall}. This motivates us to consider the correlations between multiple items in orders. 

\subsection*{Reverse-Exclude Algorithm}
This heuristic iteratively excludes products from the ground set, prioritizing those having the lowest order frequency. This approach ensures that all remaining orders can be satisfied by the selected assortment. The Reverse-Exclude algorithm follows:
\begin{algorithm}
\caption{Reverse-Exclude Algorithm}
\SetAlgoNoLine
Initialization: Eliminated SKUs: $\setE \leftarrow \emptyset$\;
\While{$|\N \setminus \setE| > K$}{
  Compute order influence  $f_i = \sum_{o\in \setO, i \in S_o} D_o $ for each $i \in \N \setminus \setE$\;
    Sort SKUs by increasing $f_i$: $L \leftarrow \text{sort}(\N \setminus \setE, \text{key}=f_i)$. Select batch size $n$ \;
    \For{$i \in L[1:n]$}{
        $\setE \leftarrow \setE \cup \{i\}$; $D_o = 0$ for $o$ such that $ i\in S_o$, where $S_o$ is the set of products in order type $o$;
    }
}
\Return{$\N \setminus \setE$}
\end{algorithm}

This algorithm starts with an empty set $\setE$ for eliminated SKUs. In each iteration, it computes the order frequency $f_i$ for each remaining SKU $i$ as the sum of demands for orders containing that SKU. 
The algorithm then sorts SKUs by increasing $f_i$ and selects a batch of $n$ SKUs with the lowest order frequency. The batch size $n$ is dynamically adjusted; it is larger when it is far from the target K and is reduced to $1$ as it approaches $K$. Each selected SKU is added to $\setE$, and the demand for all orders containing it is set to $0$. This process repeats until $|\N \backslash \setE|=K$. The algorithm returns $\N \backslash \setE$ as the final assortment.

\subsection*{Hybrid Selection Algorithm}
To leverage the complementary strengths of the ML-Top-K and Reverse-Exclude algorithms, we propose a hybrid selection approach that integrates demand forecasting with correlated purchasing patterns. The algorithm systematically combines the outputs of both methods to achieve superior performance.

The algorithm begins by identifying products selected by both the ML-Top-K and Reverse-Exclude algorithms, forming an intersection set. It then extracts products unique to each method's output, ranks them based on their respective scores, and combines the selected products using a tunable hybrid ratio $r$. The final assortment includes products from the intersection and top-ranked items from the distinct sets. 
To combine these distinct product sets, we introduce a hybrid ratio $r$ that controls the relative contribution of each method. From the ML-Top-K distinct products, we select the top $\lceil Kr \rceil$ items based on predicted demand. Similarly, we select the top $\lfloor K(1-r) \rfloor$ products from the Reverse-Exclude distinct set based on order influence scores. The final hybrid assortment combines these selections with the intersection set. This approach leverages the complementary strengths of both methods while maintaining computational efficiency.
\begin{algorithm}[htbp]
\caption{Hybrid Selection Algorithm}
\SetAlgoNoLine
\KwIn{ML-Top-K assortment $\setM$, Reverse-Exclude assortment $\setR$, hybrid ratio $r$, product rank from ML-Top-K algorithm $f_i^M$ for $i\in \setM$ and product rank from Reverse-Exclude algorithm $f_i^R$ for $i\in \setR$}
\KwOut{Final hybrid assortment $\setS$}
Initialization: $\setS \leftarrow \emptyset$\;
Compute intersection: $\setI \leftarrow \setM \cap \setR$, compute distinct products: $\setD_M \leftarrow \setM \setminus \setI$, $\setD_R \leftarrow \setR \setminus \setI$\;
Sort $\setD_M$ by decreasing $f_i^M$: $L_M \leftarrow \text{sort}(\setD_M, \text{key}=f_i^M, \text{reverse}=\text{True})$\;
Sort $\setD_R$ by decreasing $f_i^R$: $L_R \leftarrow \text{sort}(\setD_R, \text{key}=f_i^R, \text{reverse}=\text{True})$\;
Select top products: $\setS \leftarrow \setI \cup L_M[1:\lceil Kr \rceil] \cup L_R[1:\lfloor K(1-r) \rfloor]$;
\Return{$\setS$}
\end{algorithm}
The hybrid ratio $r$ serves as a tunable parameter that can be adjusted based on the relative performance of each method in different scenarios. When $r$ approaches 1, the algorithm favors ML-based demand predictions, while values closer to 0 place more emphasis on correlated purchasing patterns. This flexibility allows the algorithm to adapt to varying business contexts and data characteristics.

\section*{End-to-End Algorithm for Inventory Allocation}
\label{sec:e2e}
In this section, we introduce our new end-to-end learning framework for the inventory allocation problem. Our approach addresses the need for robust, interpretable algorithms capable of real-time decision making across large-scale SKU inventories. 

\begin{figure}[htbp]
    \centering
    \caption{The Diagram Illustrates the Recurrent Neural Network (RNN) Architecture for the End-to-End Algorithm}
    \includegraphics[width=0.75\linewidth]{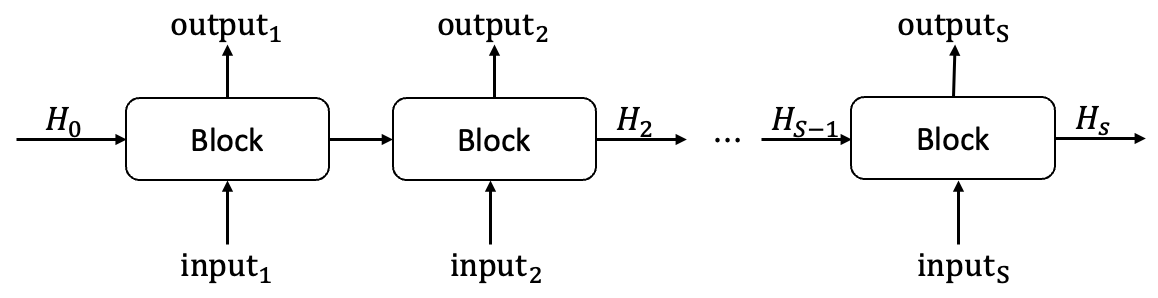}
    \label{fig:RNN}
\end{figure}

Our algorithm framework optimizes inventory allocation through a multitask, multistep decision-making procedure. Moving beyond simple single-point decisions in previous end-to-end approaches in the literature, we have implemented a process-oriented, multiperiod decision-making framework. This approach recognizes that each decision is part of a larger chain, influenced by past choices and impacting future outcomes. To capture these complex time-based relationships and utilize previous model information, we use a Recurrent Neural Network (RNN) architecture for the overall end-to-end algorithm, which we show in Figure \ref{fig:RNN}.

The framework consists of interconnected blocks, where each block is a deep learning model processing daily inputs to generate inventory management decisions. We optimize at the regional level by aggregating each RDC with its associated FDCs into a single instance, rather than solving for individual distribution centers. This regional approach enables holistic optimization that captures network-wide interdependencies while minimizing sales losses at both RDCs and FDCs and optimizing transfer costs. For each region, we calculate SS and TI based on local demand in the region, ensuring sufficient inventory to manage fluctuations during each replenishment cycle. The process begins with generating an initial SS and TI solution based on interpretable demand forecasts. We then incorporate a simulation block to minimize regional sales losses and maximize both regional and FDC-specific fulfillment rates. This simulation block helps further reduce biased inventory allocation resulting from prediction errors.

Next, we provide the details of each block.  
The inputs for each block (input$_t$) come from two sources: (1) historical sales data, marketing and promotional factors, and transfer lead times; and (2) the initial inventory state (\(H_t\)) from the previous block, including starting inventory, replenishment, and transfer amounts from RDC to FDC. The block outputs learned SS and TI levels, sales forecasts, and simulated sales for the RDC and all FDCs.

\begin{figure}[h]
    \centering   
    \caption{The Diagram Shows the Architecture in Each Block of the End-to-End Model}
    \includegraphics[width=0.95\linewidth]{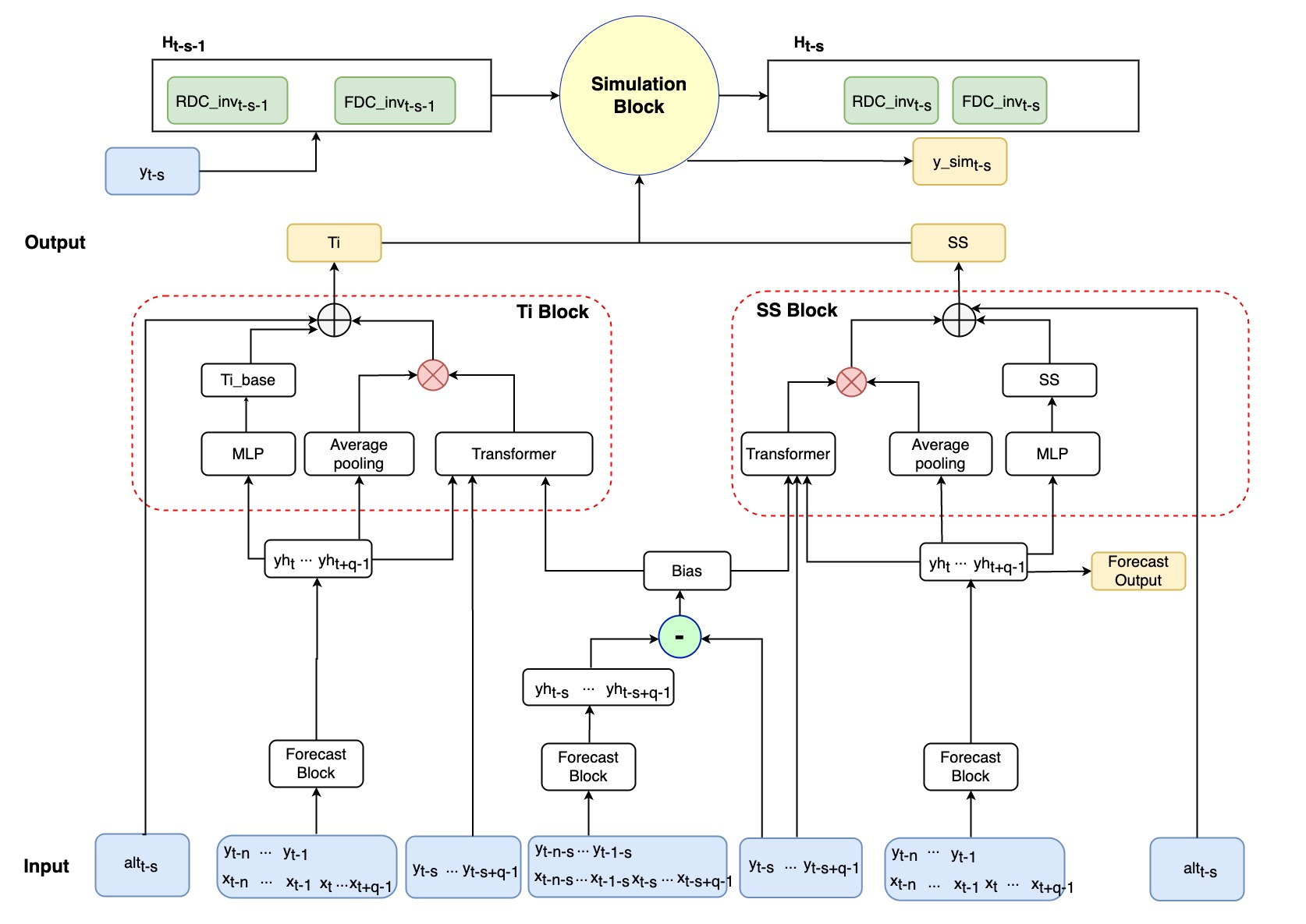}
    \label{fig:block}
\begin{minipage}{0.9\textwidth}
    \footnotesize
    Notes. $t$ denotes the time step; $n$ represents the length of historical input data; $q$ indicates the prediction horizon; and $s$ refers to the simulation duration. $x_t$ and $y_t$ denote the covariates the demand forecast at time $t$. $alt_t$ is the lead time at time $t$. $\oplus$, $\otimes$, and $\ominus$ denote different operations at each time step. $\oplus$ represents the addition operation, $\otimes$ the multiplication operation, and $\ominus$ the subtraction operation.
    \end{minipage}
\end{figure}

Figure \ref{fig:block} depicts the internal structure of each block, which has three primary components:
\begin{itemize}
\item A forecasting module that generates baseline predictions, promotional impact estimates, and adjustments based on higher-dimensional sales data (e.g., SKU-level, category-level, brand-level, and national SKU-level historical sales).
\item A TI and SS generation module that incorporates lead times, sales predictions, historical actual sales, and historical forecast accuracy to determine optimal inventory levels. This module uses MLPs and transformer structures to process the inputs and generate TI and SS values.
\item A simulation module that operates at the RDC level, allocating inventory across its network based on generated SS and TI values. The allocation strategy mirrors real-world operations by prioritizing safety stock requirements at RDCs and FDCs, followed by target inventory levels, with remaining inventory maintained at the RDC. This module simulates daily sales and inventory consumption, tracking fulfillment and lost sales. Its outputs include simulated safety stock and target inventory levels, sales forecasts for the planning horizon, and simulated daily sales across all locations.
\end{itemize}

The model training process minimizes a composite loss function $L$ consisting of three terms:
\[
L = \lambda_1 L_{\text{op}}+\lambda_2 L_{\text{sales\_pred}}+\lambda_3 L_{\text{ss}},
\]
where \(L_{\text{op}}\) represents the total operational costs from simulation including lost sales and transfer costs, $L_{\text{sales\_pred}}$ measures the prediction error between forecasted and actual sales, and $L_{\text{ss}}$
penalizes safety stock violations. Here $\lambda_1$, $\lambda_2$, and $\lambda_3$ are weighting parameters. Unlike traditional supervised learning approaches that use optimal inventory decisions as training labels, we directly optimize these operational metrics.

During the inference phase, our model generates multiday sales forecasts and recommends optimal target inventory and safety stock levels for forthcoming replenishment decisions. This approach facilitates dynamic, forward-looking inventory management that adapts to evolving demand patterns and supply chain conditions. By providing proactive recommendations, the model enables businesses to anticipate and respond to market changes efficiently, optimizing inventory levels across the distribution network.

Our framework extends previous end-to-end approaches by establishing an explicit connection between demand forecasting and inventory decisions. Although the E2E algorithm in \cite{qi2023practical} incorporates sales forecast error minimization, its architecture only implements parameter sharing between the sales forecasting and inventory prediction modules, without leveraging demand information for inventory decision making. Our framework advances this approach in two ways: it maintains parameter sharing between the sales forecasting and SS/TI generation modules while also utilizing sales forecasts directly in SS and TI predictions. This integration enables systematic analysis of model performance—when discrepancies arise between predicted SS/TI values and optimal decisions, we can analyze the corresponding sales forecasts to identify potential forecast deviations. Additionally, the simulation module provides quantitative validation through its total cost metrics. 

\section*{Numerical Experiments, Implementation, and Impact}
In this section, we present numerical results obtained from real-world data at JD.com, comparing our proposed algorithms against previously implemented methods. We also describe the implementation and impact. 
\subsection*{Numerical Experiments}
\subsubsection*{Assortment Algorithms}

We compare the performance of our proposed assortment algorithms with the previously implemented greedy algorithm based on only the order frequency from historical orders, which selects the top $K$ products with the highest historical order frequencies. 
Our experiment considers 18 FDCs, with $K$ chosen to satisfy 70\% of orders in a given historical time horizon. This fulfillment rate represents the percentage of orders where all products within the order are included in the assortment. The 70\% threshold is determined by business requirements, balancing assortment breadth with FDC capacity and inventory transfer constraints. A higher fulfillment rate would necessitate a larger number of SKUs, potentially violating these operational capacities. 

\begin{figure}
    \centering
    \caption{The Graph Illustrates Local Fulfillment Rates with Different Hybrid Ratios}
    \includegraphics[width=0.6\linewidth]{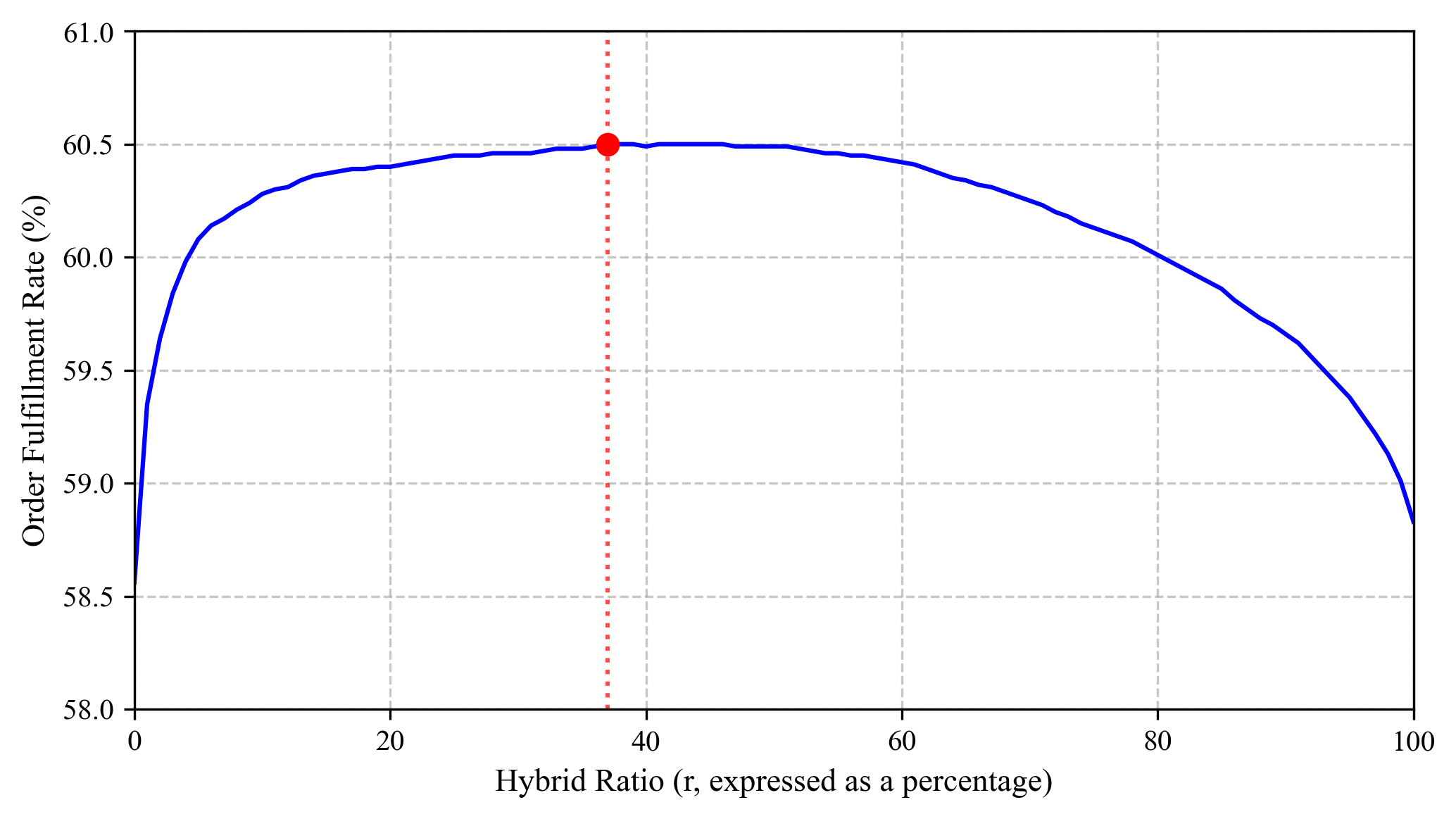}
    \label{fig:r}
\end{figure}

For training data, we use historical data from the $D$ days preceding May 1, 2024. We then generate assortments using our algorithms and the baseline method, evaluating their performance by calculating order fulfillment percentages using order data from May 1, 2024 to May 14, 2024. Our analysis includes all orders containing regular products, excluding large or special items. Note there is a hyperparameter $r$ in the hybrid selection. We select the $r$ that yields the highest local order fulfillment rate in training data. Figure \ref{fig:r} shows that the fulfillment rate initially increases rapidly at lower hybrid ratios ($r$), peaks around $r=38\%$, and then gradually decreases as $r$ further increases, suggesting diminishing returns and potential overreliance on one method at higher ratios. This suggests we should balance the solutions from the two heuristics. Thus, we set $r=38\%$ when testing the out-of-sample performance. 
Table \ref{tab:order_satisfaction} presents the fulfillment rates of the baseline method and our proposed heuristics. Our \textit{ML-Top-K} approach shows the best performance, achieving a fulfillment rate of 58.83\%. Although we cannot disclose exact figures due to data security considerations, the scale of orders fulfilled is in the range of tens of millions. This represents a significant improvement over the baseline \textit{Top-K} method, which exhibited a fulfillment rate of 58.29\%. Although the percentage difference of 0.54\% may appear modest, it translates to an additional 100,000 satisfied orders, a considerable volume in absolute terms. This improvement highlights the practical power of machine learning in inventory management. The \textit{Reverse-Exclude} algorithm also showed improvement, with a fulfillment rate of 58.56\%, surpassing the baseline by more than 50,000 orders. The advantage of this algorithm is that it is computationally efficient and does not need to train the ML model, which may be time-consuming. Furthermore, it is particularly suitable for scenarios where accurate order count predictions are challenging to obtain. The \textit{Hybrid Selection} algorithm has the highest local order fulfillment rate of 60.50\%, suggesting that leveraging both heuristics yields the best performance in practice. 
These improvements are particularly noteworthy given the scale of operations under consideration. In the context of millions of orders, even incremental percentage increases in fulfillment rates result in substantial absolute improvements. 
\begin{table}[htbp]
\centering
\caption{The Table Presents a Comparative Analysis of Local Order Fulfillment Rates Across Different Methods}
\label{tab:order_satisfaction}
\begin{tabular}{lrrr}
\hline
Method  & Fulfillment Rate (\%) \\
\hline
Top-K & 58.29 \\
ML-Top-K &  58.83 \\
Reverse-Exclude &  58.56 \\
Hybrid Selection & 60.50\\
\hline
\end{tabular}
\end{table}

\subsubsection*{Inventory Allocation}

In the experiment, we consider two RDCs, each serving six FDCs. We train the model using historical data spanning from Day $i-F-1$ to Day $i-1$, where $F$ represents the length of the historical data window used for training. We focus on three key performance metrics: FDC demand fulfillment rate (FDC Ful.), regional lost sales (Reg. Loss), and loss ratio. The loss ratio is defined as regional sales loss divided by the sales fulfilled by FDCs. This measure is crucial because transferring inventory from an RDC to an FDC can result in sales loss within the region, given the limited region capacity. The loss ratio indicates the efficiency of order fulfillment and the extent of sales loss incurred. 

We compare our proposed new end-to-end algorithm against two existing online algorithms previously used by JD.com:
\begin{enumerate}
\item Parameter Search Through Simulation: This approach simulates various inventory allocation policies at the RDC level, optimizing decisions based on demand forecasts and predefined business constraints. 
\item Linear Program (LP): This method solves a linear program for each SKU based on demand forecasts, safety stock constraints, and other business requirements to determine inventory decisions.
\end{enumerate}
Table \ref{tab:inventory_allocation} compares the performance of the three algorithms. 
The new end-to-end algorithm demonstrates superior performance in FDC demand fulfillment rate, achieving 58.59\% compared with 57.55\% for Parameter Search and 57.12\% for LP Optimization. This represents a significant improvement of 1.05 percentage points. Although our approach shows a slight increase in Regional Loss (2.86\% compared with 2.78\% for LP Optimization), the Loss Ratio remains competitive at 17.80\%, matching the Parameter Search method. The substantial improvement in the FDC fulfillment rate suggests that our method can significantly enhance local order fulfillment capabilities, which is critical for improving customer satisfaction. This improvement particularly benefits the coverage of the `211' program. 
\begin{table}[htbp]
\centering
\caption{The Table Shows a Comparison of Two Inventory Allocation Methods that JD.com Used Previously vs. the New End-to-End Method}
\label{tab:inventory_allocation}
\begin{tabular}{l|ccc}
\hline
Method & FDC Ful. (\%)& Reg. Loss (\%)  & Loss Ratio (\%)  \\
\hline
Parameter Search & 57.55 & 2.81 & 17.80 \\
LP Optimization & 57.12 & 2.78 & 17.73 \\
New End-to-End & 58.59 & 2.86 & 17.80 \\
\hline
\end{tabular}
\end{table}

\subsection*{Implementation and Impact}\label{sec:impact}
Combining the algorithms for assortment and inventory optimization, we have developed a decision support system. This system has been deployed across JD.com's inventory network in China, covering all eight RDCs, with each RDC serving 5-13 FDCs. The benefits of this system are significant: Inventory holding costs and capital utilization costs of FDCs have been reduced by approximately \$6.13 million per year. Additionally, JD.com has saved approximately \$22.32 million in inventory transfer costs. Stock availability has also improved by 0.85\%. Furthermore, the local order fulfillment rate has increased by 2.19\%, and the proportion of the speedy `211 program' has risen by 1.44\%, enhancing timeliness and customer experiences for approximately 18.61 million orders annually. These achievements show the impact of algorithmic capabilities in optimizing operations and resource allocation in the inventory network.

\section*{Conclusion}
This paper addresses the challenges of integrated assortment planning and inventory allocation in JD.com's two-level distribution network. For assortment optimization, we developed and implemented novel algorithms that effectively balance computational tractability with performance. Our proposed algorithms for assortment planning demonstrate improved local order fulfillment rates over previous methods. For inventory allocation, we propose a new end-to-end algorithm that combines forecasting, optimization, and simulation, showing superior performance in FDC demand fulfillment while maintaining competitive loss ratios.
These algorithms have been successfully deployed across JD.com's nationwide network, resulting in significant improvements in operational efficiency and customer service levels. This work shows the potential of data-driven optimization techniques in solving complex, large-scale supply chain problems. 

\bibliographystyle{informs2014}
\bibliography{papers}


\end{document}